\documentclass{IEEEtran}
\usepackage{textcomp}
\usepackage{amsmath,amsfonts,amssymb,amscd}
\usepackage[colorlinks,linkcolor=blue,citecolor=blue]{hyperref}
\usepackage[font=small]{caption}
\usepackage{latexsym}
\usepackage{graphicx,subfigure}
\usepackage{hyperref}
\usepackage{color}
\usepackage{bm}
\usepackage{cite}
\usepackage{algorithm}
\usepackage{algorithmic}
\usepackage{mathrsfs}
\usepackage{subeqnarray}
\usepackage{cases}
\usepackage{arydshln}
\usepackage{multirow}
\usepackage{geometry}
\usepackage{extarrows}

\def\BibTeX{{\rm B\kern-.05em{\sc i\kern-.025em b}\kern-.08em
		T\kern-.1667em\lower.7ex\hbox{E}\kern-.125emX}}
\newtheorem{thm}{Theorem}
\newtheorem{lem}{Lemma}
\newtheorem{remark}{Remark}
\newtheorem{ass}{Assumption}
\newtheorem{definition}{Definition}
\newtheorem{prop}{Proposition}

\newtheorem{corollary}{Corollary}

\begin{document}
	
\title{Finite-Time Convergent Algorithms for Time-Varying Distributed Optimization}

\author{Xinli~Shi, Guanghui~Wen, and Xinghuo~Yu 
	\thanks{X.~Shi is with the School of Cyber Science and Engineering, Southeast University, Nanjing 210096, China (e-mail: xinli\_shi@seu.edu.cn).}
	\thanks{G. Wen is with the the Jiangsu Provincial Key Laboratory of Networked Collective Intelligence, School of Mathematics, Southeast University, Nanjing 210096, China (e-mail: wenguanghui@gmail.com).}
    \thanks{X. Yu is with the School of Engineering, RMIT University, Melbourne, VIC 3001, Australia (e-mail: x.yu@rmit.edu.au).}
}


\markboth{Manuscript For Review}
{\;}

\maketitle

\begin{abstract}
	This paper focuses on finite-time (FT) convergent distributed algorithms for solving time-varying (TV) distributed optimization (TVDO). The objective is to minimize the sum of local TV cost functions subject to the possible TV constraints by the coordination of multiple agents in finite time. Specifically, two classes of TVDO are investigated included unconstrained distributed consensus optimization and distributed optimal resource allocation problems (DORAP) with both TV cost functions and coupled equation constraints. For the previous one, based on non-smooth analysis, a continuous-time distributed discontinuous dynamics with FT convergence is proposed based on an extended zero-gradient-sum method with a local auxiliary subsystem. Then, an FT convergent distributed dynamics is further obtained for TV-DORAP by dual transformation. Particularly, the inversion of the cost functions' Hessians is not required in the dual variables' dynamics, while another local optimization needs to be solved to obtain the primal variable at each time instant. Finally, two numerical examples are conducted to verify the proposed algorithms.
\end{abstract}
\begin{IEEEkeywords}
	Finite-time convergence, time-varying distributed optimization, resource allocation, discontinuous dynamics.
\end{IEEEkeywords}

\section{Introduction}\label{shosd:sec:introduction}

Recently, time-varying distributed optimization (TVDO) has received increasing attention as it is more practical than time-invariant distributed optimization (TIDO) in a dynamic environment where the objective function or constraints can change over the time. TVDO has found applications in many areas, such as power system, machine learning and robotics, see \cite{Simonetto2020} and references therein. In TVDO, the optimal solution can be TV and hence the traditional algorithms designed for TIDO to approach a static optimizer can not be applied directly. To solve TVDO, many discrete-time algorithms (DTAs) have been provided for solving TVDO \cite{Hosseini2016,Simonetto2017,Simonetto2018, Dixit2020,Yuan2020}. For example, prediction-correction methods are used in \cite{Simonetto2017,Simonetto2018} for solving TVDO with sample period, and the tracking error bound is related to the size of sample period. More DTAs can be found in the survey \cite{Li2022survey} and literature therein. However, in general, it is difficult for CTAs to track the TV optimal trajectory asymptotically due to the sample period, step size or errors in the local optimization. 

To track the optimal solution trajectory with vanishing errors, several continuous-time algorithms (CTAs) have been proposed for solving TV (distributed) optimization \cite{Simonetto2020,Rahili2017,Huang2020,Ning2019,Sun2017}. 
Compared with DTAs, typical CTAs can achieve fast convergence performance, such as finite-time (FT) convergence to the optima, and further inspire the discrete-time counterpart design. Besides, CTAs are also preferred for optimizing the swarm tracking behavior of a multi-robotic system with physical dynamics \cite{Simonetto2020}\cite{Rahili2017}. Moreover, CTAs can track the TV optimal solution with vanishing errors. Since CTAs exhibit such good tracking performance for TV problems, we design a CTA, with the hope that even when discretized it will show improved performance. Thus, in this work, we mainly focus on CTAs for solving TVDO, including TV distributed consensus optimization (TVDCO) and optimal resource allocation problems (TV-DORAP). For TVDCO, in \cite{Rahili2017}, both single- and double-integrator dynamics are investigated for distributed consensus optimization with TV cost functions, which can achieve asymptotic convergence to the optimal trajectory under some strict assumptions such as identical Hessians. Furthermore, TVDCO is considered for a class of nonlinear multi-agent system in \cite{Huang2020}. In \cite{Ning2019}, an edge-based distributed protocol is provided for solving TVDCO subject to identical Hessians. In \cite{Sun2017}, gradient-based searching methods are used to track the TV solution of TVDCO with quadratic objective functions.

Different from the distributed consensus optimization aiming to achieve multi-agent consensus on the optimal solution, in DORAP with coupled constraints, the optimal resource allocation of agents can be heterogeneous but should satisfy the total demand while minimizing the sum of local cost functions. For time-invariant DORAP, extensive CTAs have been proposed to approach the static optimal solution with agents' coordination \cite{Zhu2019,Jia2021TAC}. As variations always exist in the generations of renewable resources and the noncontrollable load demand in a dynamic power system \cite{Simonetto2020}, TV cost functions and constraints are more practical in the real-time DORAP. In \cite{Cherukuri2016}, a robust distributed algorithm is designed for economic dispatch and then the TV loads are discussed, under which the power mismatch is only shown to be ultimately bounded. In \cite{Wang2020TCNS}, the prediction–correction method combined with the discontinuous consensus protocol is used for DORAP with TV quadratic cost functions of identical Hessians, and then distributed average tracking (DAT) estimator-based methods are used to deal with the case with nonidentical constant Hessians. Moreover, Wang et al. provide two dynamics for handling DORAP with quadratic cost functions subject to identical and nonidentical Hessians in \cite{Wang2021} and \cite{Wang2022SMC}, respectively, as well as TV resources.

Although some CTAs have been proposed for solving TVDO including TVDCO and DORAP, there exist several issues required to be further addressed. One is that most CTAs can only achieve an asymptotic convergence to the optimal trajectory, which means that tracking errors tend to zero as the time goes to infinity. The second is that almost all the CTAs with simple structures require that the Hessians of all local cost functions are identical, and for DORAP, the cost functions are of the quadratic form, which limit the applications of the existing methods. To deal with nonidentical Hessians, a non-smooth estimator based on DAT methods can be introduced to track the global information, but it suffers from high gains and more computation/communication costs. To tackle these two issues, the task of this note is to design FT convergent CTAs for solving TVDO without identical Hessians. That is to say, the provided dynamics will track the optimal trajectory in finite time without mismatch errors. In literature, the existing distributed FT algorithms are designed for time-invariant unconstrained optimization \cite{song2016finite}. For example, an FT convergent algorithm is proposed in \cite{song2016finite} based on the continuous-time zero-gradient-sum (ZGS) method \cite{Lu2012TAC}, where the initial states should be the minimizers of the local cost functions. Recently, an FT convergent primal-dual method has been proposed in \cite{Shi2022cyber} and further applied to constrained distributed optimization including DORAP. In \cite{Hu2018Neuro}, a DAT estimator for tracking global information is used to design FT dynamics for TVDO with TV cost functions. 


As pointed out in \cite{Simonetto2020}, the best algorithm for the invariant case may be the 
worst one in the dynamic case. Similarly, the existing FT CTAs designed for TIDO might not be applied to TVDO directly, and there are few results on the FT convergent algorithms with simple structures for TVDO. In this note, we will provide two FT distributed algorithms for solving TVDO such as TVDCO and DORAP, respectively. The contributions are listed as follows.  
\begin{enumerate}
	\item For a class of TVDCO with strongly convex and smooth cost functions, a discontinuous distributed dynamics with FT convergence is proposed, based on an extended ZGS approach. Specifically, the introduced auxiliary dynamics can drive all the local states to an invariant set in finite time, on which the sum of local gradients is zero. After that, TVDCO will be solved once the multi-agent system reaches consensus in finite time by using non-smooth consensus protocol. Different from existing ZGS method based algorithms (e.g., \cite{song2016finite,Lu2012TAC}) for TIDO, the agents' initial states are not required to minimize the local functions (see \cite{song2016finite,Lu2012TAC}). Moreover, compared with exiting methods \cite{Rahili2017,Huang2020,Ning2019,Hu2018Neuro}, the provided algorithm has a simpler structure without estimating the global information and can be used for TVDO with nonidentical Hessians. However, similar to \cite{Rahili2017,Huang2020,Ning2019,Hu2018Neuro}, the inversion of Hessian is required in the proposed method.
	\item An FT convergent distributed dynamics is further obtained for DORAP with TV cost functions and coupled constraints by dual transformation. Compared with existing works \cite{Wang2020TCNS,Wang2021,Wang2022SMC}, the cost functions can be of non-quadratic form with nonidentical Hessians. Moreover, the inversion of the cost functions' Hessians is not required in the dual variables' dynamics. However, another local optimization needs to be solved to obtain the primal variable at each time instant.
	\item In the proposed CTAs, only binary information is required and the neighboring agents only need to know whether their relative position is positive or negative, which benefits the online implementation in a coarse sensing scenario \cite{Jiang2017,Chellapandi2023,Jafarian2015}. Despite such a coarse information, the FT convergence to the exact TV optimal trajectory is guaranteed.
\end{enumerate}

The rest of this paper is organized as follows. Section \ref{section-2} introduces some preliminary notations and concepts. In Section \ref{Problem}, the problem statements and main methodologies are provided. Finally, two numerical examples on TVDCO and TV DORAP are presented in Section \ref{Numberical} to verify the proposed algorithms. Conclusions are drawn in Section \ref{conclusion}.

\section{Preliminaries}\label{section-2}
\subsection{Notation and Network Representation}
$\mathbb{R}_+^n$ denotes the set of $n$-dimensional vectors with non-negative entries and $I_n\in \mathbb{R}^{n\times n}$ is an identity matrix. Let $\langle n\rangle=\{1,2,\cdots,n\}$. For $x=[x_{1},\ldots,x_{n}]^{T}\in \mathbb{R}^n$, we write the $p$-norm of $x$ as $\|x\|_p=(\sum_{i=1}^n |x_i|^p)^{\frac{1}{p}}$ and $\text{sgn}^{\alpha}(x)=[\text{sign}(x_1)|x_1|^{\alpha},...,\text{sign}(x_n)|x_n|^{\alpha}]^T$ for $\alpha\in  \mathbb{R}_+$. By default, we denote $\|x\|$ as the Euclidean norm of $x$. For the matrix $M\in \mathbb{R}^{n\times n}$, $\mathcal{N}(M)$ denotes the null space of $M$ and $\Pi_{\mathcal{N}(M)}(x)$ denotes the projection of $x$ onto $\mathcal{N}(M)$. When $M$ is positive semidefinite, $\lambda_2(M)$ denotes its smallest positive eigenvalue. 


A continuously differentiable function $f(x): \mathbb{R}^n \rightarrow \mathbb{R}$ is called 
$\mu$-strongly convex if for any $x,y \in \mathbb{R}^n$, $f(y)-f(x) -\nabla f(x)^T(y-x) \geq \frac{\mu}{2}\|y-x\|^2$.
Moreover, $f$ is said to be $\theta$-smooth if for any $x,y \in \mathbb{R}^n$, $\|\nabla f(y)-\nabla f(y)\| \leq \theta\|y-x\|$.
For a convex function  $f(x): \mathbb{R}^n \rightarrow \mathbb{R}$, its Legendre–Fenchel conjugate $f^*$ is defined by 
\begin{align}
	f^*(y) = \sup_{x \in \mathbb{R}^n}\{\langle y,x\rangle -f(x) \}, \ \forall y \in \mathbb{R}^n .
\end{align}
The following result can be found in \cite[Proposition 12.60]{Rockafellar2009}.
\begin{lem}\label{lem-fconvex}
	For a proper, lower semicontinuous and convex function $f(x): \mathbb{R}^n \rightarrow \mathbb{R}\cup \{\pm \infty\}$, $f^*$ is $\frac{1}{\sigma}$-strongly convex iff $f$ is differentiable and $\sigma$-smooth. 
\end{lem}

A weighted undirected graph is represented by $\mathcal{G}(\mathcal{V},\mathcal{E},A)$, where $\mathcal{V}=\langle N \rangle$ and $\mathcal{E}\subseteq \mathcal{V} \times \mathcal{V}$ denote the node set and edge set, respectively, and the matrix $A=[a_{ij}]_{N\times N}$ is the weighted adjacency matrix with $a_{ij}> 0$ iff $(j,i)\in \mathcal{E}$ and $a_{ii}=0, \forall i\in \langle N \rangle$. Let $\mathcal{N}_i$ denote the set of neighbors of agent $i$. 

\subsection{Finite-Time Stability}
Consider the following differential autonomous system
\begin{align}\label{auto_dynamic}
	\dot{x}(t)=f(x(t)), x(t_0)=x_0,
\end{align}
with $x\in \mathbb{R}^n$. When the right-hand side function $f: \mathbb{R}^n \rightarrow \mathbb{R}^n$ is discontinuous at some points, its Filippov solution will be investigated, which is an absolutely continuous map $x:I \rightarrow \mathbb{R}^n $ defined on an interval $I \subset \mathbb{R}$ satisfying the differential inclusion: 
\begin{align}\label{filippov}
	\dot{x}(t) \in \mathcal{F}[f](x(t))
\end{align}
with $\mathcal{F}[f](\cdot)$ being the Filippov set-vauled map \cite{ShiTAC2023}. As known, the Filippov solution to \eqref{auto_dynamic} always exists if $f$ is measurable and locally essentially bounded. The \textit{set-valued Lie derivative} of a locally Lipschitz continuous map $V: \mathbb{R}^n \rightarrow \mathbb{R}$ associated with $\mathcal{F}[f]$ at $x$ is defined as
\begin{align*}
	\widetilde{\mathcal{L}}_{\mathcal{F}[f]}V(x) \triangleq \{a \in \mathbb{R}| \exists \zeta \in \mathcal{F}[f](x) \Rightarrow \zeta^T \nu=a, \forall \nu \in  \partial V(x)\},
\end{align*}
in which $\partial V(x)$ represents the generalized gradient of $V$ \cite{ShiTAC2023}. 

The FT/FxT stability of the system \eqref{auto_dynamic} is given in Definition \ref{def-FTC}. Lemma \ref{lem-non} is helpful for the non-smooth analysis when dealing with Filippov solutions to \eqref{auto_dynamic}.
\begin{definition}\label{def-FTC}
	The system (\ref{auto_dynamic}) is \textit{FT stable} at the origin if it is Lyapunov
	stable and there exists a \textit{finite settling time} $T(x_0)$ such that
	\begin{align*}
		\lim_{t \rightarrow T(x_0)}x(t; x_0) =0 \ \text{and} \ x(t)=0, \forall t\geq T(x_0)
	\end{align*}
	for any solution $x(t; x_0)$ starting from $x_0 \in \mathcal{X}_0$. When $\mathcal{X}_0= \mathbb{R}^n$, the system (\ref{auto_dynamic}) is called globally FT stable. 
\end{definition}

\begin{lem}\cite{ShiTAC2023}\label{lem-non}
	Suppose that $V: \mathbb{R}^n \rightarrow \mathbb{R}$ is locally Lipschitz continuous and regular, and $x(t): I\subset \mathbb{R} \rightarrow \mathbb{R}^n$ is a Filippov solution satisfying \eqref{filippov}. Then, $V(x(t))$ is absolutely continuous with 
	\begin{align}
		\dot{V}(x(t)) \overset{a.e.}{\in} \widetilde{\mathcal{L}}_{\mathcal{F}[f]}V(x(t)), \ t\in I.
	\end{align}
\end{lem}


\section{Problem statement and methodology}\label{Problem}
In this section, we first provide a unified framework for designing FT convergent algorithms to solve TV centralized optimization. Then, based on an extended ZGS method, two FT convergent distributed algorithms will be designed for TVDCO and DORAP over a networked system, respectively. Meanwhile, the comparisons with existing works are discussed in remarks. 

\subsection{FT Convergent TV Optimization}
First, we focus on the following TV convex optimization
\begin{align}\label{CO-central}
	\min_{x\in \mathbb{R}^n} f_0(x,t),
\end{align}
where the objective function $f_0:\mathbb{R}^n\times \mathbb{R} \rightarrow \mathbb{R}$ is convex and further satisfies Assumption \ref{ass-f0}. 
\begin{ass}\label{ass-f0}
	The function $f_0(x,t)$ is twice continuously differentiable with invertible Hessian $H_0(x,t)$ and there exists a continuous trajectory $x^*(t)$ that solves \eqref{CO-central}.
\end{ass}

Then, we aim to provide a continuous-time dynamics $\dot{x}=u(x,t)$ such that $x(t)$ tracks $x^*(t)$ in FT $T<+\infty$, i.e.,
\begin{align}\label{xlim-cent}
	\lim_{t \rightarrow T} \|x(t)-x^*(t)\| =0.
\end{align}  
Then, to achieve \eqref{xlim-cent}, the proposed dynamics is designed as
\begin{subequations}\label{prot-cent}
	\begin{numcases}{}
		\dot{x} =-H_0^{-1}(x,t)(\varphi(z)+\frac{\partial \nabla f(x,t)}{\partial t}), \ \ \quad \label{prot-cx}  \\
		\dot{z} =-\varphi(z), \label{prot-cz}
	\end{numcases}
\end{subequations}
with $z(0)=\nabla f_0(x(0),0)$, where the function $\varphi(z)$ is chosen such that the origin is FT/FxT stable for the subsystem \eqref{prot-cz}. Several typical functions can be chosen, e.g., 
\begin{align*}
	\varphi(z) &= a\cdot z/\|z\|_r^p + b\cdot z\|z\|_r^{q-1} \ \text{or} \ a \cdot \text{sgn}^{1-p}(z) + b \cdot \text{sgn}^{q}(z)
\end{align*}
with $p \in (0,1],r\geq 1,q>1$, for which the FT/FxT stability of \eqref{prot-cz} can be shown with the Lyapunov function $V(z)=\frac{1}{2}z^Tz$ based on \cite[Lemma 1]{polyakov2011nonlinear}. Note that the subsystem \eqref{prot-cz} is only related to $z(t)$ and thus it can be freely designed with prescribed performance. In this work, we consider a class of FT stable subsystem \eqref{prot-cz}. Then, we have the following result. 
\begin{thm}\label{thm-cent}
	Suppose that Assumption \ref{ass-f0} holds and \eqref{prot-cz} is FT stable at origin with settling time $T_{\max} $. Then, with \eqref{prot-cent}, $x(t)$ will track $x^*(t)$ in the same FT $T_{\max}$. 
\end{thm}
\begin{proof}
 See Appendix A. 
\end{proof}

\begin{remark}
	In fact, from the proof of Theorem \ref{thm-cent}, \eqref{prot-cent} is equivalent to the following dynamics 
	\begin{align}\label{prot-c0}
		\dot{x} =-H_0^{-1}(x,t)(\varphi(\nabla f(x,t))+\frac{\partial \nabla f(x,t)}{\partial t}).
	\end{align}
	The introduced variable $z(t)$ is used to estimate $\nabla f(x,t)$ when the initialization coincides.
	The auxiliary system \eqref{prot-cz} will be helpful for designing FT distributed algorithms to solve TVDO in the next part. 
\end{remark}

\subsection{FT Algorithm for Solving TVDCO}

Consider an interactive network $\mathcal{G}(\mathcal{V}, \mathcal{E}, A)$ consisting of $N$ agents, each of which is equipped with a local TV cost function. Then, the multi-agent system aims to solve the following TVDO by coordinating with neighbors. 
\begin{align}\label{CO}
	\min_{x\in \mathbb{R}^n} \sum_{i=1}^Nf_i(x,t),  
\end{align}
where $f_i: \mathbb{R}^n\times\mathbb{R}_+  \rightarrow \mathbb{R}$ is the local TV cost function only known by the agent $i$. Moreover, $f_i$ is supposed to satisfy Assumptions \ref{ass-f} and \ref{ass-dft}, where $\nabla f_i(x,t)$ is the gradient of $f_i$ w.r.t. $x$. Under Assumption \ref{ass-f}, it can be deduced that there exists a unique continuous trajectory $x^*(t)$ that solves TVDO \eqref{CO}. Besides, Assumption \ref{ass-dft} covers a group of functions such as $f_i(x,t) = f_i^0(x)+g(t)x$ when $\dot{g}(t)$ is globally bounded. For example, in smart grid, $f_i(x,t)$ can be the cost of a generator with $x$ and $-g(t)$ representing the generation power and the TV electricity price, respectively. The objective of this work is to design a distributed CTA with only local information that tracks $x^*(t)$ in finite time. 

\begin{ass}\label{ass-f}
	For each local cost function $f_{i}, i\in \langle N \rangle $, $f_i$ is twice continuously differentiable w.r.t. $x$, $\underline{\theta}_i$-strongly convex and $\bar{\theta}_i$-smooth for some positive scalars $\underline{\theta}_i, \bar{\theta}_i$, i.e., $\underline{\theta}_i \bm{I}_n \leq H_i(x,t) \leq \bar{\theta}_i \bm{I}_n$, where $H_i(x,t)$ is the Hessian matrix of $f_i$. 
\end{ass}

\begin{ass}\label{ass-dft}
	For each $i \in \langle N \rangle$, $\frac{\partial \nabla f_i(x,t)}{\partial t}$ is measurable w.r.t $t$ and bounded by $\|\frac{\partial \nabla f_i(x,t)}{\partial t}\| \leq \kappa$ for some $\kappa >0$. 
\end{ass}

\begin{ass}\label{ass-G}
	The interactive graph $\mathcal{G}$ is undirected and connected.
\end{ass}

Let $x_i$ be the local copy of the system decision $x$ at agent $i$. Then, under Assumption \ref{ass-G}, \eqref{CO} can be reformulated into the following equivalent distributed optimization over the network $\mathcal{G}(\mathcal{V}, \mathcal{E})$.
\begin{align}\label{DO}
	\min_{x_i \in \mathbb{R}^n} \sum_{i=1}^Nf_i(x_i,t), \  \text{s.t.} \  x_i=x_j, \forall (i,j) \in \mathcal{E}.
\end{align}
Then, the purpose of this subsection is to design distributed dynamics which guarantees that there exists a finite time $T<\infty$ such that 
\begin{align}\label{FT}
	\lim_{t\rightarrow T} \|x_i(t) -x^*(t)\| = 0, \forall i \in \langle N \rangle.
\end{align}
In order to track the global solution $x^*(t)$ in finite time, we propose the following distributed discontinuous dynamics as an extension of \eqref{prot-cent}
\begin{subequations}\label{prot-u2}
	\small 
	\begin{numcases}{}
		\dot{x}_i =-H_i^{-1}(x_i,t)(\varphi_i(z_i)+\frac{\partial \nabla f_i(x_i,t)}{\partial t}+\alpha \phi_i(x_i, x_{\mathcal{N}_i})),  \label{prot-ux2} \\
		\dot{z}_i =-\varphi_i(z_i), \ i\in \langle N \rangle \label{prot-uz2}
	\end{numcases}
\end{subequations}
with $z_i(0) = \nabla f_i(x_i(0),0)$, where $\phi_i(x_i, x_{\mathcal{N}_i}) =\sum_{j=1}^N a_{ij}\text{sgn}(x_i-x_j) $ with $x_{\mathcal{N}_i}\triangleq [x_j]_{j\in \mathcal{N}_i}$. Note that under Assumptions \ref{ass-f}, \ref{ass-dft} and \ref{ass-z}, the right-hand side of \eqref{prot-u2} is bounded over the time and thus its Filippov solution always exists. We call \eqref{prot-u2} an extended ZGS method as the initial state $x_i(0)$ is not required to be the local minimizer of local cost function to be distinct from \cite{Lu2012TAC}. However, with the introduced auxiliary subsystem \eqref{prot-uz2}, the local states $x_i(t)$ will be driven to an invariant set where the sum of local TV gradients is zero. Particularly, when the solution to \eqref{prot-uz2} can be explicitly given as $z_i(t)= \Phi_i(t,z_0)$, $\varphi_i(z_i)$ in \eqref{prot-ux2} can be simply replaced by $\varphi_i(\Phi_i(t,z_0))$ and then \eqref{prot-uz2} can be removed. 
\begin{remark}
	For \eqref{prot-u2}, it might be restrictive that each agent knows the closed form of $\partial \nabla f_i(x_i,t)/\partial t$. However, in some applications, the local time-varying cost function is available for each agent, such as a motion control with an optimization objective, or in smart grid where smart devices should coordinate with each other to maximize the overall utility function with a TV electricity price whose rate is known beforehand \cite{Rahili2017}. The term $\partial \nabla f_i(x_i,t)/\partial t$ here is used for prediction to enhance the exact convergence to $x^*(t)$, which has been used in \cite{Rahili2017,Huang2020,Simonetto2018,Wang2020TCNS,Wang2022SMC}. As pointed in \cite{Simonetto2020}, without the predictor, a tracking error exists, possibly depending on the variation of the gradient with the time and the control gain. Besides, the discontinuous term $\phi_i$ is used for achieving consensus among agents as well as constraining $\partial \nabla f_i(x_i,t)/\partial t$. Sometimes, the noise in the channels can affect the performance of the algorithms \cite{Jiang2017,Chellapandi2023}. For $\phi_i$, as only the binary information that represents whether the neighboring agents' relative position is positive or negative is used, the presence of small noise in a coarse sensing channel may not affect the performance of algorithms when $\text{sign}(x_i-x_j)$ (i.e., $\{\pm1, 0\}$) remains the same. Detailed properties of the binary protocols can be found in \cite{Jafarian2015}.
\end{remark}

\begin{ass}\label{ass-z}
	The subsystem \eqref{prot-uz2} is FT stable at the origin, $\forall i\in \langle N \rangle$.
\end{ass}

In the next, we will show that $\sum_{i=1}^N z_i(t)$ as an estimation of $\sum_{i=1}^N \nabla f_i(x_i,t)$ will converge to zero in finite time and then all $x_i(t)$ will reach consensus on $x^*(t)$ in finite time due to $\phi_i$ in \eqref{prot-ux2}. 

Given a weighted adjacency matrix $A=[a_{ij}]_{N\times N}$ of the undirected graph $\mathcal{G}(\mathcal{V},\mathcal{E})$ with $m=|\mathcal{E}|$, we define a weighted incidence matrix $B_0=[b_{ik}]_{N \times m}$: $b_{ik} =-b_{jk}=a_{ij}$ if $i<j$ for any edge $e_k=\{i,j\} \in \mathcal{E}$. Let $\bm{x} =[x_i]_{i\in \langle N \rangle}$, $\delta =[\delta_{ij}]_{(i,j)\in \mathcal{E}}$, $\bm{z} =[z_i]_{i\in \langle N \rangle}$, $B=B_0 \otimes I_n$, $f(\bm{x},t) = \sum_{i=1}^{N}f_i(x_i,t)$, $H=\text{blkdiag}(H_1(x_i,t), \cdots, H_N(x_N,t))$ and $\varphi(\bm{z})=[\varphi_i(z_i)]_{i\in \langle N \rangle}$. Then, \eqref{prot-u2} can be rewritten in a compact form: 
\begin{subequations}\label{prot-compact}\small
	\begin{numcases}{}
		\dot{\bm{x}} =-H^{-1}(\bm{x},t)(\varphi(\bm{z})+\frac{\partial \nabla f(\bm{x},t)}{\partial t}+\alpha B\text{sgn}(B^T\bm{x})),  \label{prot-ux3}  \\
		\dot{\bm{z}} =-\varphi(\bm{z}), \label{prot-uz3}
	\end{numcases}
\end{subequations}
with $\bm{z}(0) = \nabla f(\bm{x}(0),0) = [\nabla f_i(x_i(0),0)]_{i\in \langle N \rangle}$. Then, based on the non-smooth analysis, the FT convergence of \eqref{prot-compact} is established in the following result. 
\begin{thm}\label{thm1}
	Let Assumptions \ref{ass-f}-\ref{ass-z} hold, $\underline{\theta}=\min_{i\in \langle N \rangle}\{\underline{\theta}_i\}$ and $\bar{\theta}=\max_{i\in \langle N \rangle}\{\bar{\theta}_i\}$. Under the distributed protocol \eqref{prot-u2} with $\alpha > \kappa \sqrt{\frac{N\bar{\theta}}{\underline{\theta}\lambda_2(B^TB)}}$, $z_i(t)$ will converge to 0 in finite time and then $x_i(t)$ will track the optimal trajectory $x^*(t)$ to TVDO \eqref{DO} in finite time. 
\end{thm}
\begin{proof}
	See Appendix B. 
\end{proof}


Actually, by the Eq. \eqref{V-relax}, the strong convexity assumption and Assumption \ref{ass-dft} can be further relaxed as Assumption \ref{ass-relax}, which is more general than \cite[Ass. 3.8]{Rahili2017}, \cite[Ass. 4]{Huang2020}, \cite[Ass. 2]{Ning2019} and \cite{Simonetto2018}. It can be easily verified that when $f_i$ is strongly convex and $\|\partial \nabla f_i(x,t)/\partial t\|$ is bounded, Assumption \ref{ass-relax} holds directly. With the alternative Assumption \ref{ass-relax}, we have the following result. 
\begin{ass}\label{ass-relax}
	For any $i,j \in \langle N \rangle$, there exists $\varpi >0$ such that $\|H_i^{-1}(x_i,t)\frac{\partial \nabla f_i(x_i,t)}{\partial t}-H_j^{-1}(x_j,t)\frac{\partial \nabla f_j(x_j,t)}{\partial t}\|_1 \leq \varpi$. 
\end{ass}
\begin{prop}\label{prop-relax}
	Let Assumptions \ref{ass-G}-\ref{ass-relax} hold. Suppose that $f_i(x_i,t)$ has invertible Hessian and is $\theta_i$-smooth. Let $\bar{\theta}=\max_{i\in \langle N \rangle}\{\bar{\theta}_i\}$ and $\bar{a}= \max_{(i,j) \in \mathcal{E}}\{a_{ij}\}$. Under the distributed protocol \eqref{prot-u2} with $\alpha > \frac{m\varpi\bar{a}\bar{\theta}}{\lambda_2(B^TB)}$, $z_i(t)$ will converge to 0 in finite time and then $x_i(t)$ will track the optimal trajectory $x^*(t)$ to TVDO \eqref{DO} in finite time. 
\end{prop}
\begin{proof}
See Appendix C. 
\end{proof}
\begin{remark}In \cite{Simonetto2018}, a dual prediction–correction DTA is developed for solving equation constrained TV optimization and then its distributed implementation is given by using constraint matrix (i.e., the
	incidence/Laplacian matrix) to ensure the local consensus. In \cite{Rahili2017}, an asymptotically convergent CTAs is provided for solving \eqref{DO} subject to TV cost functions with identical Hessians, based on an FT DAT dynamics to track the global information. The work \cite{Rahili2017} first makes all the states achieve consensus in finite time and then the TV optimization is solved by a centralized Newton–Raphson (NR) dynamics, i.e., $\dot{x}=-(\nabla^2f(x))^{-1}\nabla f(x)$, which we call ``consensus+central NR" dynamics. Differently, non-identical Hessians are allowable in the proposed ZGS-based algorithm \eqref{prot-u2} and the state consensus term is inside the local NR dynamics, which we call ``consensus within local NR" dynamics. The advantage of the ZGS-based algorithm lies in that the TVDCO is solved once the consensus of local states is reached as the ZGS manifold $\mathcal{M}_0$ is invariant. By the convergence analysis, one can only obtain an upper bound of the finite settling time applicable for all considered TVDCO in worst case as most results for FT algorithms do. For specific problems, such as DAT system, the derived control gain condition in Theorem \ref{thm1} is less conservative than that (i.e., $\alpha > \kappa \sqrt{2N/\lambda_2(B^TB)}$) obtained in \cite[Lemma 1]{Wang2020TCNS}, for which $\bar{\theta} = \underline{\theta}=1$. To improve the convergence rate, one can design the auxiliary system \eqref{prot-uz2} to achieve the expected finite settling time when the local states reach the $\mathcal{M}_0$, and then enlarge the control gain $\alpha$ to further reduce the finite consensus time. However, large gain might result in high chattering magnitude with discretization.
\end{remark}

\begin{remark}
	Assumption \ref{ass-dft} used in Theorem \ref{thm1} is more general than existing works \cite{Rahili2017,Ning2019,Huang2020}, and holds for an important class of functions such as $f_i(x_i,t) = (a_ix_i+b_i(t))^2/2$, where $a_i$ is heterogeneous among agents and only the boundedness of $\|\dot{b}_i(t)\|$ is required, e.g, $b_i(t)=t$ and $\text{sin}(t)$. However, in \cite{Rahili2017}, for the distributed method with neighbors' position, with identical Hessians, i.e., $a_i=a_j$, the boundedness of $\|b_i(t)-b_j(t)\|$, $\|\dot{b}_i(t)-\dot{b}_j(t)\|$ and $\|\ddot{b}_i(t)-\ddot{b}_j(t)\|$ is required to satisfy \cite[Assumption 3.8]{Rahili2017}. Similarly, identical Hessians are used in \cite{Ning2019,Huang2020}. For DAT-based method provided in \cite{Rahili2017}, to deal with non-identical Hessians, the boundedness of $\|\frac{d \nabla f_i(x_i,t)}{dt}\|$, $\|\frac{d H_i(x_i,t)}{dt}\|$ and $\|\frac{d (\partial f_i(x_i,t)/\partial t)}{dt}\|$ is required and known for the FT convergence of the DAT estimators. In fact, these items not only depend on $\dot{b}_i(t)$ and $\ddot{b}_i(t)$, but also rely on the evolution of local states $x_i$, and thus are difficult to estimate. For a general class of quadratic cost functions, e.g., $f_i(x_i,t) = \frac{1}{2}x_i^TH_i(t)x_i+b_i^T(t)x_i$, as $\nabla f_i(x_i,t) = H_i(t)x_i+b_i(t)$ also depends on the value of $x_i$, the value $\kappa$ in Assumption \ref{ass-dft} is not easy to be known generally. However, if $x_i(t)$ is known to be upper bounded, with the boundedness of $\|\dot{H}_i(t)\|$ and $\dot{d}_i(t)$, Assumption \ref{ass-dft} holds directly. 
\end{remark}

\subsection{FT Algorithm for Solving TV-DORAP}
In this subsection, the TV-DORAP will be studied over a multi-agent system, where each agent is associated with a TV local cost function. All the involved agents will coordinate with each other to minimized the whole cost function subject to a coupled equation constraint. The DORAP can be formulated as follows:
\begin{align}\label{DORAP}
	\min \sum_{i=1}^N f_i(x_i,t) \quad \text{s.t.}\quad  \sum_{i=1}^{N}x_i=d(t), 
\end{align}  
where $x_i \in \mathbb{R}^n$ is the local state of agent $i$, $d(t) \in \mathbb{R}^n$ is the TV total resource that will be allocated among all agents, and $f_i: \mathbb{R}^n\times \mathbb{R}_+\rightarrow \mathbb{R}$ is the TV local cost function of agent $i$. We further suppose that $d(t)= \sum_{i=1}^{N}d_i(t)$, where $d_i(t)\in \mathbb{R}^n$ can be accessed only by the $i$th agent. Note that \eqref{DORAP} is different from \eqref{DO} where all the local states are required to achieve consensus on the optimal solution. 

For the TV-DORAP \eqref{DORAP}, by introducing the Lagrangian multiplier associated with the equality constraint, the TV Lagrangian function can be written as 
\begin{align}
	\mathcal{L}(\bm{x},\lambda,t) = \sum_{i=1}^N f_i(x_i,t) +\lambda \sum_{i=1}^{N}(d_i(t)-x_i),
\end{align}
where $\bm{x}=[x_i]_{i\in \langle N \rangle}\in \mathbb{R}^{Nn}$. Then, by \cite[Sec. 5]{Boyd2004}, one can derive the Lagrangian dual function $g(\lambda,t)$ 
\begin{align*}\small
	g(\lambda,t) = \inf_{\bm{x}}	\mathcal{L}(\bm{x},\lambda,t)
	= \sum_{i=1}^N (\lambda d_i(t) -f_i^*(\lambda,t)) =\sum_{i=1}^N g_i(\lambda,t),
\end{align*} 
where $g_i(\lambda,t)= \lambda d_i(t) -f_i^*(\lambda,t)$ can be regarded as the local distributed dual function. 
Under Lemma \ref{lem-fconvex} and Assumption \ref{ass-f}, it can be concluded that $g_i$ is $\frac{1}{\bar{\theta}_i}$-strongly concave and $\frac{1}{\underline{\theta}_i}$-smooth. 
Moreover, as $\nabla f_i^*$ is an inverse map of $\nabla f_i$ and $\nabla^2f_i(x_i(\lambda), t)$ is positive definite, it holds that $\nabla^2 f_i^*(\lambda,t) = (\nabla^2f(x_i(\lambda,t),t))^{-1}$, where $x_i(\lambda,t)=\arg \sup_{x_i}(\lambda x_i - f_i(x_i,t))$. Besides, $\nabla g_i(\lambda,t) = d_i(t) -x_i(\lambda,t)$. 
Then, the dual problem of \eqref{DORAP} is given by 
\begin{align}\label{DO-dual}
	\max_{\lambda \in \mathbb{R}^n} \sum_{i=1}^N g_i(\lambda,t), 
\end{align}
which can be transferred to the TVDCO \eqref{DO} as
\begin{align}\label{DO-dual2}
	\min_{\bm{\lambda} \in \mathbb{R}^n} \sum_{i=1}^N -g_i(\lambda_i,t), \  \text{s.t.} \  \lambda_i=\lambda_j, \forall (i,j) \in \mathcal{E}
\end{align}
with $\bm{\lambda} = [\lambda_i]_{i\in \langle N \rangle}\in \mathbb{R}^{Nn}$.
Then, the distributed algorithm \eqref{prot-u2} can be used to solve \eqref{DO-dual2} as follows. 
\begin{subequations}\label{prot-dual}
	\small 
	\begin{numcases}{}
		\dot{\lambda}_i =-H_i(x_i,t)(\varphi_i(z_i)-\frac{\partial \nabla g_i(\lambda_i,t)}{\partial t}
		+\alpha\phi_i(\lambda_i, \lambda_{\mathcal{N}_i})),		\label{prot2-ulam2}  \\
		\dot{z}_i =- \varphi_i(z_i),  \label{prot2-uz2} \\
		x_i(\lambda_i,t) = \arg \sup_{x_i \in \mathbb{R}^n}(\lambda_i x_i - f_i(x_i,t)), \  \forall i\in \langle N \rangle \label{prot2-ux2} 
	\end{numcases}
\end{subequations}
with $z_i(0) = -\nabla g_i(\lambda_i(0),0)=x_i(\lambda_i,0)-d_i(0)$ and $\frac{\partial \nabla g_i(\lambda_i,t)}{\partial t} = \dot{d}_i(t)-\frac{\partial x_i(\lambda_i,t)}{\partial t}$. Take $f_i(x_i,t) = \frac{1}{2}a_ix_i^2 + b_i(t)x_i+c_i(t)$ as an example with $x_i \in \mathbb{R}$. One can derive that $x_i(\lambda_i,t) = \frac{\lambda_i-b_i(t)}{a_i}$ and $\frac{\partial \nabla g_i(\lambda_i,t)}{\partial t} = \dot{d}_i(t)+\frac{\dot{b}_i(t)}{a_i}$. To ensure the boundedness of $\frac{\partial \nabla g_i(\lambda_i,t)}{\partial t}$, Assumption \ref{ass-dt} is imposed on $d_i(t)$. Then, one can show the boundedness on $\frac{\partial \nabla g_i(\lambda_i,t)}{\partial t}$ in Lemma \ref{lem-dbound} under Assumptions \ref{ass-f}, \ref{ass-dft} and \ref{ass-dt}. 
\begin{ass}\label{ass-dt}
	There exists $\delta >0$ such that $\|\dot{d}_i(t)\| \leq \delta, \forall i \in \langle N \rangle$. 
\end{ass}
\begin{lem}\label{lem-dbound}
	Let Assumptions \ref{ass-f}, \ref{ass-dft} and \ref{ass-dt} hold. It satisfies that 
	\begin{align}
		\|\frac{\partial \nabla g_i(\lambda_i,t)}{\partial t}\| \leq \frac{\kappa}{\underline{\theta}_i}  +\delta, \ \forall i \in \langle N \rangle. 
	\end{align}
\end{lem}
\begin{proof}
	See Appendix D. 
\end{proof}

By applying Theorem \ref{thm1} and Lemma \ref{lem-dbound}, one can derive the following corollary. 
\begin{corollary}\label{cor-1}
	Suppose that Assumptions \ref{ass-f}-\ref{ass-dt} hold. Let $\underline{\theta}=\min_{i\in \langle N \rangle}\{\underline{\theta}_i\}$ and $\bar{\theta}=\max_{i\in \langle N \rangle}\{\bar{\theta}_i\}$. Consider the distributed protocol \eqref{prot-dual} with $\alpha > (\frac{\kappa}{\underline{\theta}}  +\delta) \sqrt{\frac{N\bar{\theta}}{\underline{\theta}\lambda_2(B^TB)}}$, $z_i(t)$ will converge to 0 in finite (resp. fixed) time and then $x_i(t)$ will track the optimal trajectory $x_i^*(t)$ to the TV DORAP \eqref{DORAP} in finite time. 
\end{corollary}

\begin{remark}
	In \cite{Wang2020TCNS} and \cite{Wang2021}, several asymptotically convergent algorithms are provided for solving DORAP with TV quadratic cost functions. Moreover, TV demand $d_i(t)$ is considered in \cite{Wang2022SMC}, where both $\|d_i(t)\|_{\infty}$ and $\|\dot{d}_i(t)\|_{\infty}$ are required to be bounded. In \cite{Wang2020TCNS} and \cite{Wang2022SMC}, an estimator based on DAT method is first provided for tracking the global information in finite time, and then a distributed CTA is designed to track the TV optimal solution asymptotically. Differently, the proposed algorithm \eqref{prot-dual} in this note has a simple form and can track the optimal trajectory in finite time. Moreover, compared with existing works \cite{Wang2020TCNS,Wang2021,Wang2022SMC}, the cost functions can be of non-quadratic form with non-identical Hessians and the inverses of Hessian of $f_i$ is not required in \eqref{prot-dual} from a dual perspective, which extends the existing results. However, another local optimization needs to be solved to obtain the primal variable at each time instant.
\end{remark}

\begin{figure*}[t]
	\centering
	\begin{minipage}[t]{0.22 \textwidth}
		\centering
		{\includegraphics[width=0.87 \textwidth]{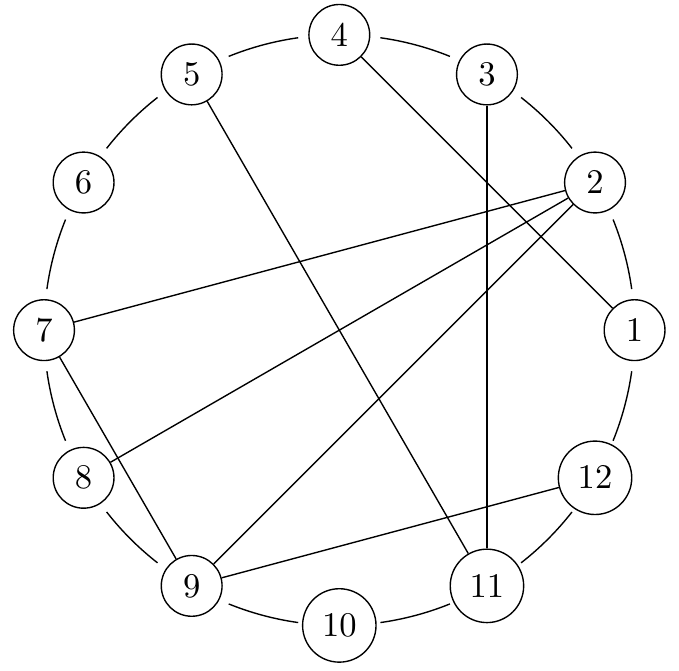}}
		\caption{\scriptsize{The network with 12 agents.}}
		\label{Fig1}
	\end{minipage}
	\begin{minipage}[t]{0.252 \textwidth}
		\centering
		{\includegraphics[width=0.95\textwidth]{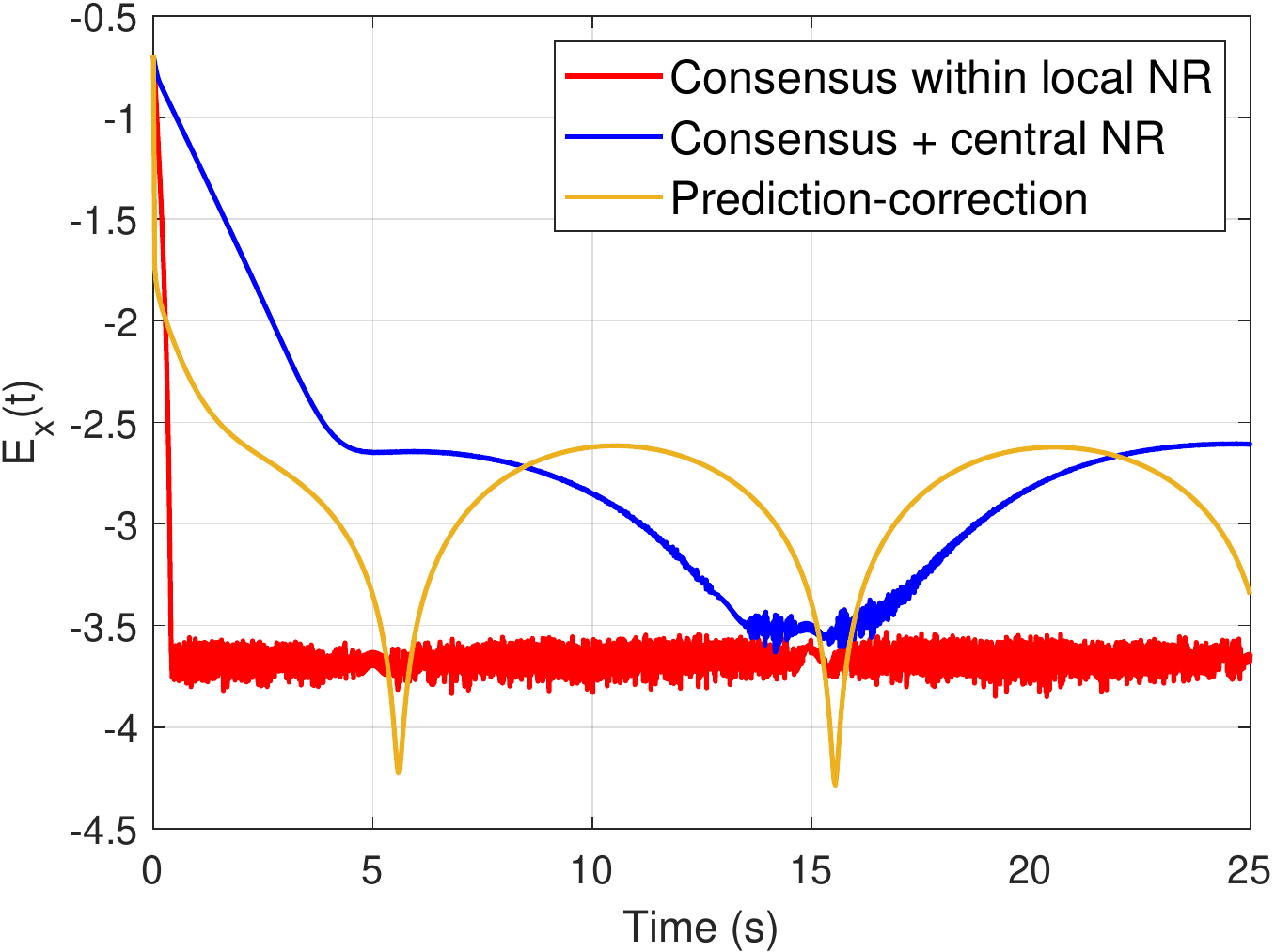}}
		\caption{\scriptsize Case 1: Values of $E_x(t)$.}
		\label{case1}
	\end{minipage}
	\begin{minipage}[t]{0.254 \textwidth}
		\centering
		{\includegraphics[width=0.93\textwidth]{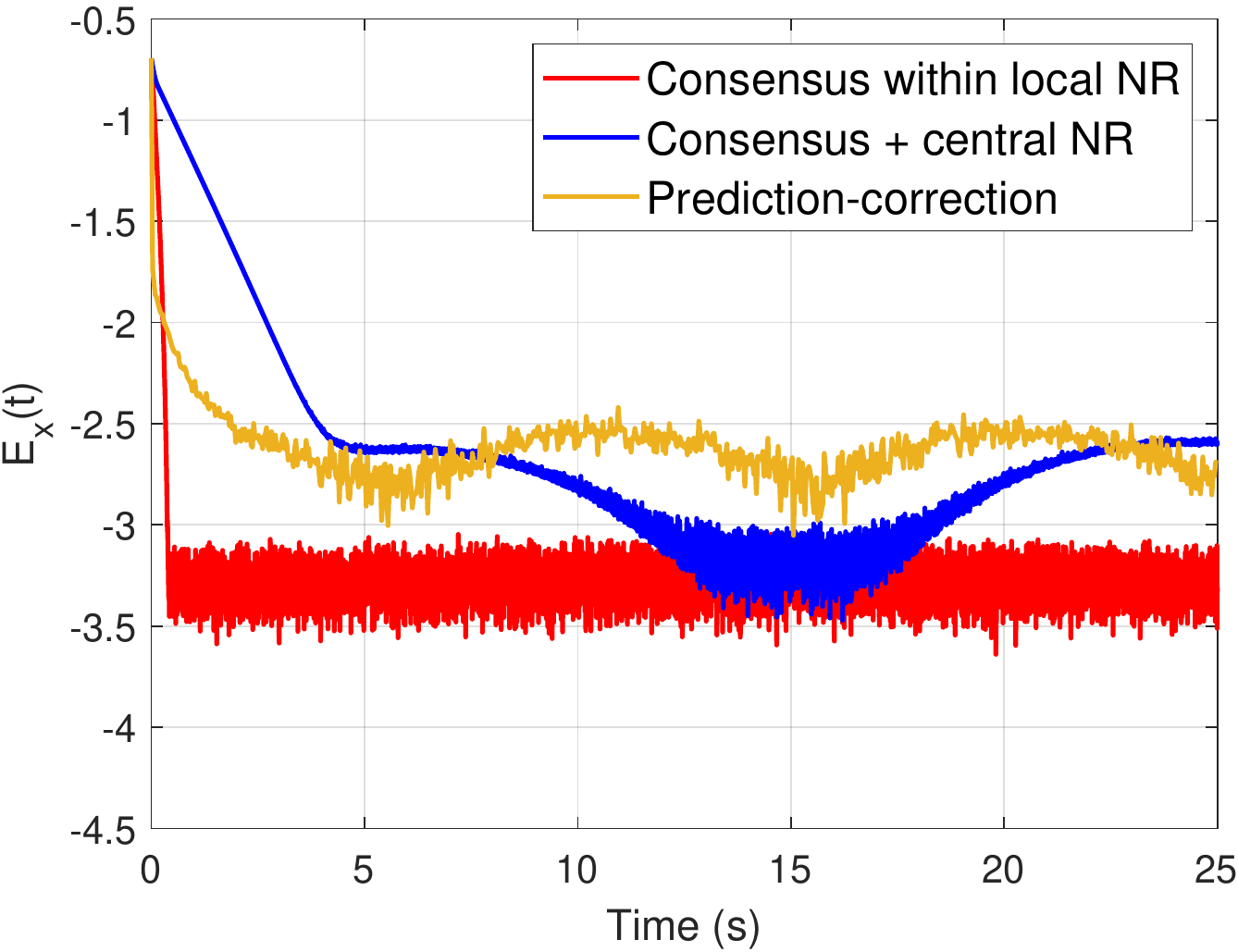}}
		\caption{\scriptsize Case 1: $E_x(t)$ with noise.}
		\label{case1-noise}
	\end{minipage} 
	\begin{minipage}[t]{0.25\textwidth}
		\centering
		{\includegraphics[width=0.94 \textwidth]{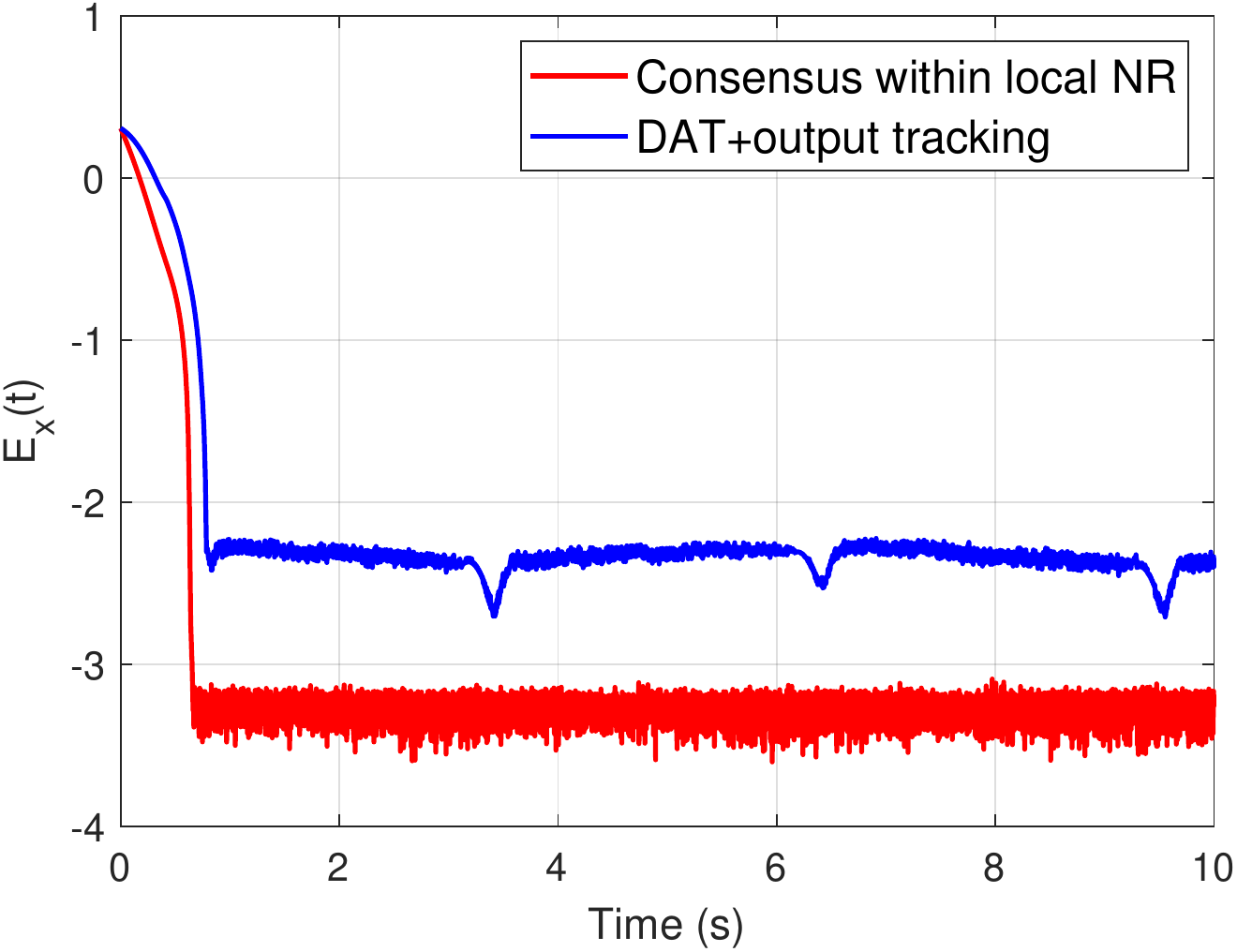}}
		\caption{\scriptsize Case 2: Values of $E_x(t)$.}
		\label{Fig-case2}
	\end{minipage}
\end{figure*}

\section{Numerical examples}\label{Numberical}
In this section, two case studies will be conducted to demonstrate the effectiveness of the proposed algorithms \eqref{prot-u2} and \eqref{prot-dual}, respectively. We consider a multi-agent system consisting of $N=12$ agents over the communication network $\mathcal{G}(\langle N \rangle,\mathcal{E})$ as shown in Fig. \ref{Fig1}, for which the edge weights are all ones. Then, one can calculate that $\lambda_2(B^TB)=1.239$ with $B$ being the incidence matrix of $\mathcal{G}$. All the experiments are conducted using Matlab R2021b on a 2.9 GHz Intel Core i7.

For the first case study, we consider the following TV binary logistic regression model with $l_2$-regularization:
\begin{align}\label{LR}
	\min_{x} f(x) = \sum_{i=1}^N(\log(1+\exp(-l_iy_i^T(t)x)) + \frac{\beta_i}{2} \|x\|_2^2),
\end{align} 
where $l_i \in \{-1,1\}$, $y_i(t) \in \mathbb{R}^2$ and $\beta_i$ is the regularization parameter, an integer randomly chosen in $[1\ 8]$. In the simulation, the TV sample data $y_i(t)=(1+\text{sin}(wt))y_i(0)$ with $w = \pi/10$ and $y_i(0)$ chosen randomly in $[0 \ 1]^2$. The proposed ``consensus within local NR" dynamics \eqref{prot-u2} with $\alpha=4$ and $\varphi(z_i)=10\text{sgn}^{0.5}(z_i)$ is compared with the ``consensus+central NR" dynamics \cite[Alg. (9)]{Rahili2017} and ``prediction–correction" algorithm \cite{Simonetto2018} with Laplacian constraint matrix. The performance of three algorithms is testified online in a real scenario. The former two CTAs are implemented by Euler discretization with step $h = 0.4 $ ms as the computation/communication of each update consumes less than $0.4$ ms. The ``prediction–correction" algorithm is executed with $P=3$ prediction steps and $C=3$ correction steps, for which the sampling period is set as $h=0.03$ s considering the computation/communication time for 3 prediction and correction operations using the \texttt{fmincon} function in MATLAB. Note that two rounds of communication are required for each prediction/correction operation. The primal states $x_i(0), i\in \mathcal{V}$ are randomly chosen in $[0 \ 1]^2$. Let $E_x(t) = \log_{10}(\frac{1}{N}\sum_{i=1}^{N}\|x_i(t)-x^*(t)\|)$ denote the tracking error between real-time states and the TV optimal solution $x^*(t)$. The simulation results are shown in Fig. \ref{case1}. It can be seen that the proposed algorithm converges to a small neighborhood of $x^*(t)$ in finite time less than 1s and then stay therein, which is faster and more robust than the other two algorithms. However, the final chattering exists due to the discontinuous controller and Euler discretizaion. For comparison, the states generated by \cite[Alg. (9)]{Rahili2017} asymptotically converge to a larger neighborhood of $x^*(t)$ compared with \eqref{prot-u2} partially due to the fact that the local Hessians are not identical \cite{Rahili2017}. For the ``prediction–correction" algorithm, the local solutions also fluctuate around $x^*(t)$ possibly due to the sample period and the TV optimal solution. To testify the effects of noisy and inexact information on the performance of algorithms, the Gaussian noise $\epsilon \sim \mathcal{N}(0, 10^{-4})$ is added to each link and $\partial \nabla f_i(x_i,t)/\partial t$. The simulation results shown in Fig. \ref{case1-noise} indicate that the accuracy of all algorithms is reduced because of the noise. However, due to the binary protocol, the convergence behavior of the former two algorithms is less affected than that of the last one.

For the second case study, the DORAP \eqref{DORAP} is considered. For comparison, as the existing algorithms mostly consider DORAP with quadratic cost functions, here the cost functions $f_i$ are given by $f_i(x_i,t) = 0.5(2+0.1i) x_i^2 +\sin(0.1it)x_i+\sin(0.6it) $. The local TV demand is set as $d_i(t) = i + \text{sin}(t+i\pi/N), \forall i\in \langle N \rangle$. The proposed ``consensus within local NR" dynamics \eqref{prot-dual} with $\alpha=5$ and $\varphi(z_i)=10\text{sgn}^{0.5}(z_i)$ is compared with \cite[Alg. (6)]{Wang2022SMC} (with all gains equal to 3), named by ``DAT+output tracking" method as two DAT dynamics are firstly performed to track the global information and then an output tracking dynamics is given for ensuring the KKT condition. The algorithms are implemented by Euler discretization with step $h = 0.2 $ ms by considering the computation/communication time in each update. Let $E_x(t) = \log_{10}(\frac{1}{N}\sum_{i=1}^{N}\|x_i(t)-x_i^*(t)\|)$ where $x_i^*(t)$ is the optimal TV solution. The simulation results shown in Fig. \ref{Fig-case2} indicate that the states driven by the proposed algorithm converge to a small neighborhood of $x^*(t)$ in finite time less than 1s, with a higher accuracy than \cite[Alg. (6)]{Wang2022SMC}. Significantly, in this case (only quadratic functions are considered in \cite{Wang2022SMC}), both the computation/communication and storage cost of the proposed algorithm at each instant are less than those of \cite[Alg. (6)]{Wang2022SMC} as more auxiliary variables are used and communicated in \cite[Alg. (7)]{Wang2022SMC}.

\section{Conclusion}\label{conclusion}
In this paper, both TV centralized and distributed optimization have been investigated. For TV centralized optimization, a unified approach is given for designing FT convergent algorithms. Then, a distributed discontinuous dynamics with FT convergence is proposed for solving TVDCO based on an extended ZGS method with the local auxiliary dynamics. 
Furthermore, the proposed distributed dynamics is used for solving TV-DORAP with additional TV coupled equation constraints by dual transformation. Significantly, the cost functions can be of non-quadratic form with non-identical Hessians and the inversion of the cost functions' Hessians is not required in the dual variables' dynamics. However, another local optimization needs to be solved to obtain the primal variable at each time instant. Further work will focus on TVDO with general TV constraints over switching networks, as well as the convergence analysis on the discretized version of the proposed algorithm.

\appendices 

\section*{Appendix}

\subsection{Proof of Theorem \ref{thm-cent}}\label{App-C0}
	Considering the dynamics \eqref{prot-cent}, one can derive that 
\begin{align}
	d \nabla f(x,t)/d t = H_0(x,t)\dot{x} + \partial \nabla f(x,t)/\partial t = -\varphi(z) = \dot{z}.
\end{align}
Since $z(0)= \nabla f(x(0),0)$, $z(t)$ and $ \nabla f(x,t)$ coincide, i.e., $z(t) \equiv \nabla f(x,t), \forall t \geq 0$. As \eqref{prot-cz} is FT stable at origin with settling time $T_{\max} $, it indicates that $\nabla f(x(t),t)$ will converge to zero in finite time $T_{\max}$, i.e., \eqref{xlim-cent} holds with $T= T_{\max}$.

\subsection{Proof of Theorem \ref{thm1}}\label{App-CA}
By Assumption \ref{ass-z}, there exists a finite time $T_1$ such that $\bm{z}(t)=0, \forall t\geq T_1$. Based on \eqref{prot-u2}, it can be derived that 
\begin{align*}
	\frac{d}{dt} \sum_{i=1}^{N}\nabla f_i(x_i,t) &= \sum_{i=1}^{N}(H_i(x_i,t)\dot{x}_i + \frac{\partial \nabla f_i(x_i,t)}{\partial t} ) \\
	&= -\sum_{i=1}^{N} \varphi_i(z_i) = \sum_{i=1}^{N} \dot{z}.
\end{align*}
Since $\sum_{i=1}^{N}z_i(0) = \sum_{i=1}^{N}\nabla f_i(x_i(0),0)$, then it holds that 
\begin{align}
	\sum_{i=1}^{N} z_i(t) = \sum_{i=1}^N \nabla f_i(x_i(t),t), \quad \forall t\geq 0. 
\end{align} 
From the previous analysis, $\bm{z}(t)=0$ for any $t\geq T_1$, which indicates that $\bm{x}(t) \in \mathcal{M}_0 \triangleq \{\bm{x}(t): \sum_{i=1}^N \nabla f_i(x_i(t),t) =0\}$. Then, to show that $x_i(t)$ tracks $x^*(t)$ in finite time, it is sufficient to prove that all the local states $x_i$ achieve consensus in finite time. By \eqref{prot-uz3}, for $t\geq T_1$, one can restrain that $\varphi(\bm{z})=0$ since $\bm{z}=0$. Then, \eqref{prot-ux3} reduces to 
\begin{align*}
	\dot{\bm{x}} =-H^{-1}(\bm{x},t)(\frac{\partial \nabla f(\bm{x},t)}{\partial t}+\alpha B\text{sgn}(B^T\bm{x})).
\end{align*} 
Denote $\hat{\bm{x}}=B^T\bm{x} $, $g(t) = \frac{\partial \nabla f(\bm{x},t)}{\partial t}$ and $H(t) = H(\bm{x},t)$. It holds that 
\begin{align}\label{x-deri0}
	\dot{\hat{\bm{x}}} =-B^T H^{-1}(t)(g(t)+\alpha B\text{sgn}(\hat{\bm{x}})).
\end{align} 
Denote the right-hand side of \eqref{x-deri0} as $X_r$. By the sum rule, we have 
\begin{align}
	\mathcal{F}[X_r](\bm{x},t) \subseteq -B^T H^{-1}(t)(\mathcal{F}[g](t)+\alpha B\mathcal{F}[\text{sgn}](\hat{\bm{x}})).
\end{align}
Consider the Lyapunov function $V_2(t) = \|\hat{\bm{x}}\|_1$. By Lemma \ref{lem-non}, it holds that $\dot{V}_2 \overset{a.e.}{\in} \widetilde{\mathcal{L}}_{\mathcal{F}[X_r]}V_2$, i.e., there exist $\gamma \in \mathcal{F}[g](t)$ and $\zeta \in \mathcal{F}[\text{sgn}](\hat{\bm{x}})$ such that 
\begin{align}
	\dot{V}_1 \overset{a.e.}{=} -\nu^TB^T H^{-1}(t)(\gamma+\alpha B \zeta_2), \ \forall \nu \in \partial \|\hat{\bm{x}}\|_1. 
\end{align}
Note that $\partial \|\hat{\bm{x}}\|_1 = \mathcal{F}[\text{sgn}](\bm{\hat{x}}(t))$. Then, by choosing $\nu = \zeta$, it gives
\begin{align}
	\dot{V}_1 & \overset{a.e.}{=} -\zeta^TB^T H^{-1}(t)(\gamma+\alpha B \zeta) \nonumber \\
	& \overset{a.e.}{=} -\zeta^TB^T H^{-1}(t)\gamma-\alpha \zeta^TB^TH^{-1}(t) B\zeta  \label{V-relax}\\
	& \leq  -(\alpha - \frac{\varepsilon}{2})\zeta^TB^T H^{-1}(t) B\zeta+\frac{1}{2\varepsilon}\gamma^TH^{-1}(t) \gamma \nonumber\\
	& \leq  -(\alpha - \frac{\varepsilon}{2})\zeta^TB^T H^{-1}(t) B\zeta+\frac{N\kappa^2}{2\varepsilon \underline{\theta}},\nonumber
\end{align}
where the first inequality is obtained by applying $a^Tb \leq \frac{1}{2}(\varepsilon a^Ta+\frac{1}{\varepsilon}b^Tb)$ with $a = -H^{-\frac{1}{2}}(t) B\zeta$ and $b =H^{-\frac{1}{2}}(t)g(t)$, and the second inequality is due to Assumptions \ref{ass-f} and \ref{ass-dft} by setting $\underline{\theta}=\min_{i\in \langle N \rangle}\{\underline{\theta}_i\}$. Denote $M(t)= B^T H^{-1}(t) B$. It holds that  
\begin{align}
	\zeta^T M(t) \zeta \geq \lambda_2(H^{-1}(t)) \zeta^T B^T B \zeta \geq  \frac{1}{\bar{\theta}}\zeta^T B^T B \zeta
\end{align} 
with $\bar{\theta}=\max_{i\in \langle N \rangle}\{\bar{\theta}_i\}$. As $B^T B$ is positive semidefinite, then it holds that 
\begin{align}
	\zeta^TB^T B\zeta \geq \lambda_2(B^T B) \|\zeta -\Pi_{\mathcal{N}(B)}(\zeta) \|^2
\end{align}
since $\mathcal{N}(B^T B) = \mathcal{N}(B)$. Let $\eta = \Pi_{\mathcal{N}(B)}(\zeta)$. As $\hat{\bm{x}} =B^T\bm{x}$, then $\hat{\bm{x}}^T\eta =0$. Due to $\zeta =\text{sign}(\hat{\bm{x}})$, when $\hat{\bm{x}} \neq 0$, there exists at least one entry $k$ such that $\zeta_k=\text{sign}(\hat{\bm{x}}_k)\neq 0$ and $\zeta_k \neq \text{sign}(\eta_k) $, which implies that $\|\zeta -\eta \| \geq \|\zeta_k -\eta_k \|\geq 1$ and thus $\zeta^TB^T B\zeta\geq \lambda_2(B^TB)$. Then, we have   
\begin{align*}
	\dot{V}_1  \leq  -(\alpha - \frac{\varepsilon}{2})\frac{1}{\bar{\theta}}\lambda_2(B^TB) +\frac{N\kappa^2}{2\varepsilon\underline{\theta}} \triangleq -\varrho,
\end{align*}
for which the right side is minimized at $\varepsilon^*=\kappa \sqrt{\frac{N}{\underline{\theta}\bar{\theta}\lambda_2(B^TB)}} $.
By choosing $\alpha > \kappa \sqrt{\frac{N\bar{\theta}}{\underline{\theta}\lambda_2(B^TB)}}$, one gets $\varrho>0$ and $\dot{V}_1  \leq -\varrho$, implying that $V_1(t)$ will converge to zero in finite time $T_2$. In other words, all the local states $x_i$ reach consensus in $T_2$. As $\bm{x}(t) \in \mathcal{M}_0$ for $t\geq T_1$, one can conclude that $x_i(t) = x^*(t)$ for any $t\geq T_2$.

\subsection{Proof of Proposition \ref{prop-relax}}\label{App-CB}
	Following the proof of Theorem \ref{thm1} and denoting $\bar{a}= \max_{(i,j) \in \mathcal{E}}\{a_{ij}\}$, Eq. \eqref{V-relax} can be further relaxed as
\begin{align*}
	\dot{V}_1 & = -\zeta^TB^T H^{-1}(t)\gamma-\alpha \zeta^TB^TH^{-1}(t) B\zeta  \\
	& \leq  -\alpha \zeta^TB^T H^{-1}(t) B\zeta+m\varpi\bar{a}  \\
	& \leq  - \frac{\alpha}{\bar{\theta}}\lambda_2(B^TB) +m\varpi\bar{a}.
\end{align*}
By choosing $\alpha > \frac{m\varpi\bar{a}\bar{\theta}}{\lambda_2(B^TB)}$, the conclusion can be shown similar to the proof of Theorem \eqref{thm1}.

\subsection{Proof of Lemma \ref{lem-dbound}}\label{App-CC}
	By the definition of $x_i(\lambda_i,t)$ in \eqref{prot2-ux2}, it holds that $\lambda_i-\nabla f_i(x_i(\lambda_i,t),t)=0 $, from which it can be derived that 
\begin{align}
	H_i(x_i,t) \frac{\partial x_i(\lambda_i,t)}{\partial t} +  \frac{\partial \nabla f_i(x_i,t)}{\partial t} =0. 
\end{align}  
Then, we have $\frac{\partial x_i(\lambda_i,t)}{\partial t} =- H_i^{-1}(x_i,t)\frac{\partial \nabla f_i(x_i,t)}{\partial t} $ and hence by applying Assumptions \ref{ass-f}, \ref{ass-dft} and \ref{ass-dt}, for any $i \in \langle N \rangle$, it holds that 
\begin{align*}
	\|\frac{\partial \nabla g_i(\lambda_i,t)}{\partial t}\| &= \|\frac{\partial x_i(\lambda_i,t)}{\partial t}-\dot{d}_i(t)\|\\
	& \leq \|H_i^{-1}(x_i,t)\frac{\partial \nabla f_i(x_i,t)}{\partial t}\| +\|\dot{d}_i(t)\|  \leq \frac{\kappa}{\underline{\theta}_i}  +\delta.
\end{align*}


\bibliographystyle{ieeetr}

\begin{thebibliography}{10}
	\bibitem{Simonetto2020}
	A. Simonetto, E. Dall’Anese, S. Paternain, G. Leus, and G. B. Giannakis,
	“Time-varying convex optimization: Time-structured algorithms and applications,” {\em Proc. of the IEEE}, vol. 108, no. 11, pp. 2032–2048, Nov. 2020.
	
	
	\bibitem{Hosseini2016}	
	S. Hosseini, A. Chapman, and M. Mesbahi, “Online distributed convex optimization on dynamic networks," {\em IEEE Trans. Autom. Control}, vol. 61, no. 11, pp. 3545-3550, Nov. 2016.
	
	\bibitem{Simonetto2017}	
	A. Simonetto, A. Koppel, A. Mokhtari, G. Leus, and A. Ribeiro, “Decentralized prediction-correction methods for networked time-varying convex optimization,” {\em IEEE Trans. Autom. Control}, vol. 62, no. 11, pp. 5724-5738, Nov. 2017.
	\bibitem{Simonetto2018}	
	A. Simonetto, ``Dual prediction–correction methods for linearly constrained time-varying convex programs,'' {\em IEEE Trans. Autom. Control}, vol. 64, no. 8, pp. 3355-3361, 2018.
	\bibitem{Dixit2020}	
	R. Dixit, A. S. Bedi, and K. Rajawat, ``Online learning over dynamic graphs via distributed proximal gradient algorithm,'' {\em IEEE Trans. Autom. Control}, vol. 66, no. 11, pp. 5065-5079, 2020.
	\bibitem{Yuan2020}	
	K. Yuan, W. Xu, and Q. Ling, ``Can primal methods outperform primal-dual methods in decentralized dynamic optimization,'' {\em IEEE Trans. Signal Process.}, vol. 68, pp. 4466-4480, 2020.
	\bibitem{Li2022survey}	
	X. Li, L. Xie, and N. Li, A survey of decentralized online learning. arXiv preprint arXiv:2205.00473.
	\bibitem{Simonetto2020}
	E. Dall'Anese, A. Simonetto, S. Becker, and L. Madden, ``Optimization and learning with information streams: Time-varying algorithms and applications,'' {\em IEEE Sig. Process. Mag.}, vol. 37, no. 3, pp. 71-83, 2020.
	
	
	
	\bibitem{Rahili2017}
	S.~Rahili and W.~Ren, ``Distributed continuous-time convex optimization with
	time-varying cost functions,'' {\em IEEE Trans. Autom. Control}, vol.~62, no.~4, pp.~1590--1605, Apr. 2017.
	
	
	
	\bibitem{Huang2020}
	B.~Huang, Y.~Zou, Z.~Meng, and W.~Ren, ``Distributed time-varying convex
	optimization for a class of nonlinear multiagent systems,'' {\em IEEE Trans. Autom. Control}, vol.~65, no.~2, pp.~801--808, Feb. 2020.
	
	\bibitem{Ning2019}
	B.~Ning, Q.-L. Han, and Z.~Zuo, ``Distributed optimization for multiagent
	systems: An edge-based fixed-time consensus approach,'' {\em IEEE Trans. Cybern.}, vol.~49, no.~1, pp.~122--132, Jan. 2019.
	
	\bibitem{Sun2017}
	C. Sun, M. Ye, and G. Hu, “Distributed time-varying quadratic optimization for multiple agents under undirected graphs,” {\em IEEE Trans. Autom. Control}, vol. 62, no. 7, pp. 3687–3694, Jul. 2017.
	
	
	
	
	\bibitem{Zhu2019}
	Y. Zhu, W. Ren, W. Yu, and G. Wen, “Distributed resource allocation over directed graphs via continuous-time algorithms,” {\em IEEE Trans. Syst., Man, Cybern., Syst.}, vol. 51, no. 2, pp. 1097-1106, Feb. 2021.
	
	\bibitem{Jia2021TAC}
	W. Jia, N. Liu, and S. Qin, “An adaptive continuous-time algorithm for nonsmooth convex resource allocation optimization," {\em IEEE Trans. Autom. Control}, vol. 67, no. 11, pp. 6038-6044, Nov. 2022.
	
	
	\bibitem{Cherukuri2016}
	A. Cherukuri and J. Cort\'es, “Initialization-free distributed coordination for economic dispatch under varying loads and generator commitment,” {\em Automatica}, vol. 74, pp. 183--193, 2016.
	
	\bibitem{Wang2020TCNS}
	B.~Wang, S.~Sun, and W.~Ren, ``Distributed continuous-time algorithms for
	optimal resource allocation with time-varying quadratic cost functions,''
	{\em IEEE Trans. Control of Netw. Syst.}, vol.~7, no.~4, pp.~1974--1984, Dec. 2020.
	
	\bibitem{Wang2021}
	B.~Wang, Q.~Fei, and Q.~Wu, ``Distributed time-varying resource allocation
	optimization based on finite-time consensus approach,'' {\em IEEE Contr. Syst. Lett.}, vol.~5, no.~2, pp.~599--604, Apr. 2021.
	
	\bibitem{Wang2022SMC}
	B.~Wang, S.~Sun, and W.~Ren, ``Distributed time-varying quadratic optimal
	resource allocation subject to nonidentical time-varying Hessians with
	application to multiquadrotor hose transportation,'' {\em IEEE Trans. Syst., Man, Cybern., Syst.}, vol. 52, no. 10, pp. 6109-6119, Oct. 2022. 
	
	\bibitem{Bai2018CDC}
	L. Bai, C. Sun, Z. Feng, and G. Hu, “Distributed continuous-time resource allocation with time-varying resources under quadratic cost functions,” {\em in Proc. IEEE Conf. Decis. Control}, Miami Beach, FL, USA,
	Dec. 2018, pp. 823–828.
	
	\bibitem{song2016finite}
	Y. Song and W. Chen, “Finite-time convergent distributed consensus optimisation over networks,” {\em IET Control Theory Appl.}, vol. 10, no. 11, pp. 1314-1318, 2016.
	
	
	
	
	
	
	\bibitem{Lu2012TAC}
	J.~Lu and C.~Y. Tang, ``Zero-gradient-sum algorithms for distributed convex
	optimization: The continuous-time case,'' {\em IEEE Trans. Autom. Control}, vol.~57, no.~9, pp.~2348--2354, Sept. 2012.
	
	
	\bibitem{Shi2022cyber}
	X. Shi, X. Xu, X. Yu, and J. Cao, “Finite-time convergent primal-dual gradient dynamics with applications to distributed optimization,” {\em IEEE Trans. Cybern.}, vol. 53, no. 5, pp. 3240-3252, 2023. 
	
	\bibitem{Hu2018Neuro}
	Z. Hu and J. Yang, ``Distributed finite-time optimization for second order continuous-time multiple agents systems with time-varying cost function,'' {\em Neurocomputing}, vol. 287, pp. 173--184, 2018.
	
	\bibitem{Jiang2017}	
	Z. Jiang, K. Mukherjee and S. Sarkar, ``Convergence and noise effect analysis for generalized gossip-based distributed optimization," 2017 American Control Conference (ACC), Seattle, WA, USA, 2017, pp. 4353-4358.
	\bibitem{Chellapandi2023}	
	V. P. Chellapandi, A. Upadhyay, A. Hashemi, and S. H. \.{z}ak, ``On the convergence of decentralized federated learning under imperfect information sharing," {\em IEEE Contr. Syst. Lett.}, vol. 7, pp. 2982-2987, 2023.
	
	\bibitem{Jafarian2015}	
	M. Jafarian and C. De Persis, ``Formation control using binary information," {\em Automatica}, vol. 53, pp. 125-135, 2015.
	
	\bibitem{Rockafellar2009}	
	R. T. Rockafellar and R. J.-B. Wets, {\em Variational Analysis}, Grundlehren
	der Mathematischen Wissenschaften, Springer-Verlag, 2009.
	
	
	
	\bibitem{ShiTAC2023}
	X. Shi, G. Wen, J. Cao, and X. Yu, “Finite-time distributed average tracking for multi-agent optimization with bounded inputs,” {\em IEEE Trans. Autom. Control}, vol. 68, no. 8, 2023.
	
	\bibitem{Boyd2004}
	S. P. Boyd and L. Vandenberghe. Convex optimization. Cambridge university press, 2004.
	
	\bibitem{polyakov2011nonlinear}
	A.~Polyakov, ``Nonlinear feedback design for fixed-time stabilization of linear
	control systems,'' {\em IEEE Trans. Autom. Control}, vol.~57, no.~8, pp.~2106--2110, Aug. 2011.
	
	
	
	
	
	
	
	
\end{thebibliography}

\end{document}